\documentstyle[12pt]{article}

\newtheorem{predfn}{{\bf Definition}}

\newenvironment{dfn}{\begin{predfn}{\hspace{
               -.3em}{\bf.}}}{\end{predfn}}
\newtheorem{prelem}{{\bf Lemma}}

\newenvironment{lem}{\begin{prelem}{\hspace{
               -.3em}{\bf.}}}{\end{prelem}}
\newtheorem{prethm}{{\bf Theorem}}

\newenvironment{thm}{\begin{prethm}{\hspace{
               -.3em}{\bf.}}}{\end{prethm}}
\newtheorem{precor}{{\bf Corollary}}

\newenvironment{cor}{\begin{precor}{\hspace{
               -.3em}{\bf.}}}{\end{precor}}
\newtheorem{preprop}{{\bf Proposition}}

\newenvironment{prop}{\begin{preprop}{\hspace{
               -.3em}{\bf.}}}{\end{preprop}}
\newtheorem{preproof}{{\bf Proof.}}

\newenvironment{proof}[1]{\begin{preproof}{\rm
               #1}\hfill{$\Box$}}{\end{preproof}}

\def\emline#1#2#3#4#5#6{%
       \put(#1,#2){\special{em:moveto}}%
       \put(#4,#5){\special{em:lineto}}}
\def\newpic#1{}

\newcommand{\Deg}{\mbox{\rm deg}\,}
\newcommand{\mod}{\mbox{\rm mod}\ }
\newcommand{\shc}{semi--H--cordial}
\newcommand{\hc}{H--cordial}
\newcommand{\zmc}{zero--M--cordial}
\renewcommand{\iff}{if and only if }
\newcommand{\hkc}{H$_k$--cordial}
\newcommand{\htc}{H$_2$--cordial}

\title{\LARGE\sf A Note on ``H-Cordial Graphs''}
\author{M.~Ghebleh and R. Khoeilar}
\date{Institute for Studies in Theoretical Physics\\
and Mathematics (IPM)\vspace{.2ex}\\
and\vspace{.2ex}\\
Department of Mathematical Sciences\\
Sharif University of Technology\\
Tehran, Iran
}

\begin{document}
\addtolength{\baselineskip}{1mm}
\maketitle

\begin{abstract}
The concept of an \hc\ graph is introduced by I.~Cahit in 1996
(Bulletin of the ICA). But that paper has some gaps and invalid
statements. We try to prove the statements whose proofs in Cahit's
paper have problems, and also we give counterexamples for the
wrong statements. We prove necessary and sufficient conditions for
\hc ity of complete graphs and wheels and \htc ity of wheels,
which are wrongly claimed in Cahit's paper.
\end{abstract}

\section{Introduction}

\hc\ graphs is introduced by I.~Cahit in~\cite{cahit}, and as he
claims they can be useful to construct Hadamard matrices since any
$n\times n$ Hadamard matrix gives an \hc\ labeling for the
complete bipartite graph $K_{n,n}$. But of course the inverse is
not necessarily true. Unfortunately Cahit's paper has many wrong
statements and proofs. For example the second part of
``Lemma~2.3'' obviously is not true. To see that consider trees in
Figure~\ref{exmpl}.

Or the definition of a \zmc\ labeling there is not valid, since no
such labelings exist. Here we try to recover that paper by fixing
some wrong proofs, and restating some statements. In this section
we mention some definitions and preliminaries which are referred
to throughout the paper.

We consider simple graphs (which are finite, undirected, with no loops or
multiple edges). For the necessary definitions and notation we refer the
reader to standard texts, such as~\cite{west}.

For a {\sf labeling} of a graph $G$ we mean a map $f$ which
assigns to each edge of $G$ an element of $\{-1,+1\}$. If a
labeling $f$ is given for a graph $G$, for each vertex $v$ of $G$
we define $f(v)$ to be the sum of the labels of all edges having
$v$ as an endpoint. In other words $f(v)=\sum_{e\in I(v)}f(e)$,
where $I(v)$ is the set of all edges incident to $v$. For an
integer $c$ we define $e_f(c)$ to be the number of edges having
label~$c$, and similarly $v_f(c)$ is the number of vertices having
the label~$c$. The following lemma which states a simple but
essential relation is immediate.

\begin{lem}
\label{essential}
If $f$ is an assignment of integer numbers to the edges and vertices of a
given graph $G$ such that for each vertex $v$, $f(v)=\sum_{e\in I(v)}f(e)$,
then $\sum_{v\in V(G)}f(v)=2\sum_{e\in E(G)}f(e)$.
\end{lem}

\begin{dfn}
A labeling $f$ of a graph $G$ is called {\sf\hc}, if there exists a positive
constant $K$, such that for each vertex $v$, $|f(v)|=K$, and the
following two conditions are satisfied,\\
\centerline{$|e_f(1)-e_f(-1)|\le 1$\ \ \ and\ \ \ $|v_f(K)-v_f(-K)|\le 1$.}
A graph $G$ is called to be {\sf\hc}, if it admits an \hc\ labeling.
\end{dfn}

The following lemma provides the most--used technique of the present paper.

\begin{lem}
\label{nechc}
If a graph $G$ with $n$ vertices and $m$ edges is \hc\ then $m-n$ is even.
\end{lem}

\begin{proof}{
Since $|v_f(-1)-v_f(1)|\le 1$, if $n$ is even we have
$v_f(-1)=v_f(1)$, and by Lemma~\ref{essential} we have
$\sum_{e\in E(G)}f(e)=0$. This implies that $m$ is even.
Using a similar argument one can prove that if $m$ is even then $n$
is also even.
}\end{proof}

If $G$ is a tree, $m-n=-1$, so we have the following.

\begin{cor}
\label{hctree}
No \hc\ tree exists.
\end{cor}

\section{Trees}

Now that an \hc\ tree do not exist, we can study \shc ity of trees instead,
which is a weaker condition than \hc ity.

\begin{dfn}
A labeling $f$ of a tree $T$ is called {\sf\shc}, if for each vertex $v$,
$|f(v)|\le 1$, and $|e_f(1)-e_f(-1)|\le 1$, and $|v_f(1)-v_f(-1)|\le 1$.
A tree $T$ is called to be {\sf\shc}, if it admits a \shc\ labeling.
\end{dfn}

In~\cite{cahit} ``Lemma~2.3'' states that if $T$ is a tree such
that each of its vertices has odd degree, then $n_I$, the number
of internal vertices of $T$ satisfies the following
$$n_I\equiv\left\{\begin{array}{ll} 0\ (\mod 2)&;{\rm if}\ \
n\equiv 2\ (\mod 4)\\ 1\ (\mod 2)&;{\rm if}\ \ n\equiv 0\ (\mod
4)\\
\end{array}\right.$$
We mentioned in the last section that this statement is not true.
One can see this by two simple examples.

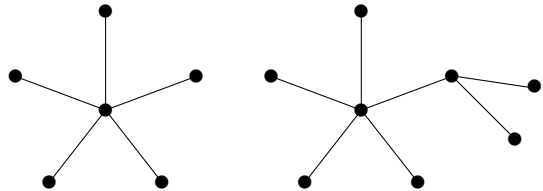
\begin{figure}[hbt]
\centering
\special{em:linewidth 0.4pt}
\unitlength 0.30mm
\linethickness{0.4pt}
\begin{picture}(235.83,81.17)
\put(43.33,78.67){\circle*{6.25}}
\put(83.33,49.67){\circle*{6.25}}
\put(68.33,2.67){\circle*{6.25}}
\put(18.33,2.67){\circle*{6.25}}
\put(3.33,49.67){\circle*{6.25}}
\put(43.33,34.67){\circle*{6.25}}
\emline{43.33}{78.67}{1}{43.33}{34.67}{2}
\emline{83.33}{49.67}{3}{43.33}{34.67}{4}
\emline{68.33}{2.67}{5}{43.33}{34.67}{6}
\emline{18.33}{2.67}{7}{43.33}{34.67}{8}
\emline{3.33}{49.67}{9}{43.33}{34.67}{10}
\put(156.66,78.67){\circle*{6.25}}
\put(196.66,49.67){\circle*{6.25}}
\put(181.66,2.67){\circle*{6.25}}
\put(131.66,2.67){\circle*{6.25}}
\put(116.66,49.67){\circle*{6.25}}
\put(156.66,34.67){\circle*{6.25}}
\emline{156.66}{78.67}{11}{156.66}{34.67}{12}
\emline{196.66}{49.67}{13}{156.66}{34.67}{14}
\emline{181.66}{2.67}{15}{156.66}{34.67}{16}
\emline{131.66}{2.67}{17}{156.66}{34.67}{18}
\emline{116.66}{49.67}{19}{156.66}{34.67}{20}
\put(224.66,22.00){\circle*{6.25}}
\put(233.33,45.34){\circle*{6.25}}
\emline{232.66}{44.67}{21}{196.66}{50.00}{22}
\emline{196.66}{50.00}{23}{224.66}{22.00}{24}
\end{picture}
\caption{Counterexamples for ``Lemma~2.3'' of Cahit's paper}
\label{exmpl}
\end{figure}

Each vertex in each of trees in Figure~\ref{exmpl} has odd degree.
The tree on the left has
six vertices and one internal vertex, and the right one has eight vertices
and two internal ones.

The proofs of Lemma~2.1, Lemma~2.2, and Theorem~2.5 in~\cite{cahit}
have serious problems, for example in some of them ``Lemma~2.3'' is used
which we showed that is not valid.
But the statements of Lemma~2.1, Lemma~2.2, and Theorem~2.5 are true
and we prove all of these in the following theorem.

\begin{thm}
\label{lbl}
A tree $T$ is \shc, \iff it has an odd number of vertices.
\end{thm}

\begin{proof}{
Suppose that $T$ has an even number of vertices,
and $f$ is a \shc\ labeling for $T$.
For each vertex $v$, we have $f(v)\in\{-1,0,1\}$, so if $\Deg v$ is even
then $f(v)=0$. Since $T$ has an even number of odd vertices
by a similar argument as in the proof of Lemma~\ref{nechc},
this means that $T$ has an even number
of edges which contradicts the hypothesis.

Now assume that $T$ has an odd number of vertices. We find a \shc\ labeling
$f$ for $T$ using an algorithm.

{\bf Algorithm.} Define two variables $S$ and $a$, where $S$
is a set and $a$ is a number. Initially we have $S=E(T)$ and $a=1$.
Update $S$ and $a$ using the following two steps while $S\not=\emptyset$.

\begin{enumerate}
\item Suppose that $e_1,e_2,\ldots,e_p$ is the longest path in $S$.
For each $1\le i\le p$ define $f(e_i)$ to be $(-1)^ia$ and then delete
$e_i$ from $S$.
\item If $\sum_{e\in E(T)\setminus S}f(e)\not=0$ then set $a$ to be equal
to it, otherwise set $a=1$.
\end{enumerate}

We claim that $f$ is a \shc\ labeling for $T$.
First note that after each execution of the two operations, we have
$a\in\{-1,1\}$, because
in the \hbox{$k$--th} execution if $p$ is even, then $f(e_1)+\ldots+f(e_p)=0$
and $\sum_{e\in E(T)\setminus S}f(e)$ do not change. Otherwise we have
$f(e_1)+\ldots+f(e_p)=-a$ and $\sum_{e\in E(T)\setminus S}f(e)$ changes
to $0$ or $-a$. So for each $e\in E(T)$ we have $f(e)\in\{-1,1\}$.
Now if the edges incident to $v$ are completely deleted from $S$ in the
\hbox{$k$--th} execution, then we have $f(v)=0$ before the \hbox{$k$--th}
execution and $|f(v)|\le 1$ after the \hbox{$k$--th} execution till the
end of algorithm. On the other hand we see that
$\sum_{v\in V(T)}f(v)=2\sum_{e\in E(T)}f(e)=0$, so $v_f(-1)=v_f(1)$.
}\end{proof}

In Lemma~2.6 of \cite{cahit} a special case of the following proposition
is stated but the proof in~\cite{cahit} has problem. We prove the
statement in a rather simple way.

\begin{prop}
Let $T$ be a tree with an even number of vertices. There exists an
edge--labeling $f$ of $T$, such that $|e_f(-1)-e_f(1)|=1$, $|f(v)|\le 1$
for each vertex $v$ in $T$, and $|v_f(-1)-v_f(1)|=2$.
\end{prop}

\begin{proof}{
Suppose that $v$ is a leaf in $T$. Add a new vertex $w$ and a new edge
$vw$ to $T$. The resulting tree has a \shc\ labeling $f$ by Theorem~\ref{lbl}
and the restriction of $f$ to $T$ is what we look for.
}\end{proof}

``Theorem~2.8''~\cite{cahit} states that a tree is \shc\ \iff it
has an even number of vertices. We have proved the opposite in
Theorem~\ref{lbl}.

\section{\hc\ graphs}

The concept of a \zmc\ labeling defined in~\cite{cahit} is useful
while one tries to find an \hc\ labeling for a given graph, There
are some wrongs on the concept occurred in~\cite{cahit}. For
example the definition of a \zmc\ labeling given in~\cite{cahit}
is not useful, because no such labelings exist! But the following
is what one expects for a \zmc\ labeling.

\begin{dfn}
A labeling $f$ of a graph $G$ is called {\sf\zmc}, if for each vertex $v$,
$f(v)=0$. A graph $G$ is called to be {\sf\zmc}, if it admits a
\zmc\ labeling.
\end{dfn}

In~\cite{cahit} the definition has an additional condition
$|e_f(-1)-e_f(1)|\ge 1$. However Lemma~\ref{essential} for a \zmc\ labeling
$f$, implies that $\sum f(e)=0$, hence  $e_f(-1)=e_f(1)$. So no graph
may have a \zmc\ labeling in sense of~\cite{cahit}.

The usefulness of the above definition appears when one tries to
find an \hc\ labeling for a given graph $G$. If $H$ is a \zmc\
subgraph of $G$, then \hc ity of $G\setminus E(H)$ simply implies
\hc ity of $G$. We will do so in the proof of Theorem~\ref{KnHC}.

It is immediate from the definition that a graph is \zmc, \iff each of
its components is a \zmc\ graph.
In the following theorem we give a characterization of connected
\zmc\ graphs.

\begin{thm}
A connected graph $G$ is \zmc\ \iff it is Eulerian and it has an even
number of edges.
\end{thm}

\begin{proof}{
Obviously each vertex in a \zmc\ graph must have even degree, and because
$e_f(1)=e_f(-1)$, it must have an even number of edges. On the other hand
if $G$ is an Eulerian graph with even number of edges, one can label edges
in order on an Eulerian tour, alternately by $+1$ and $-1$, to attain a
\zmc\ labeling.
}\end{proof}

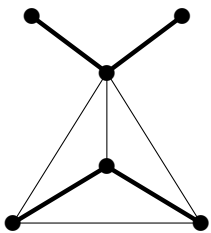
\begin{figure}[hbt]
\centering
\unitlength .50mm
\linethickness{0.4pt}
\begin{picture}(54.58,59.58)
\special{em:linewidth 0.4pt}
\put(3.33,3.33){\circle*{4.00}}
\put(53.33,3.33){\circle*{4.00}}
\put(28.33,43.33){\circle*{4.00}}
\put(28.33,18.33){\circle*{4.00}}
\put(8.33,58.33){\circle*{4.00}}
\put(48.33,58.33){\circle*{4.00}}
\emline{28.33}{43.33}{1}{3.33}{3.33}{2}
\emline{3.33}{3.33}{3}{53.33}{3.33}{4}
\emline{53.33}{3.33}{5}{28.33}{43.33}{6}
\emline{28.33}{43.33}{7}{28.33}{18.33}{8}
\special{em:linewidth 1.7pt}
\emline{3.33}{3.33}{9}{28.33}{18.33}{10}
\emline{28.33}{18.33}{11}{53.33}{3.33}{12}
\emline{8.33}{58.33}{13}{28.33}{43.33}{14}
\emline{28.33}{43.33}{15}{48.33}{58.33}{16}
\end{picture}
\caption{Counterexample for ``Theorem~3.1'' of Cahit's paper}
\label{oddm1}
\end{figure}

``Theorem~3.1''~\cite{cahit} gives a necessary condition for a labeling $f$
of a connected graph $G$, to be \hc, that is the number
of vertices labeled $-1$ must be even. We show this is not always true
by an example. The graph shown in Figure~\ref{oddm1} is our example. For an
\hc\ labeling of this graph, one can assign $-1$ to thin edges and $+1$ to
thick ones.

In ``Theorem~3.6''~\cite{cahit} it is claimed that if $n\equiv 0\ (\mod 8)$,
then the complete graph $K_n$ has an \hc\ labeling $f$ such that
$|f(v)|=n-1, \forall v\in V(G)$.
This is not true since if such $f$ exists, then all edges incident with
a specified vertex must have the same label. This implies that all edges of
$K_n$ must have the same label, which is impossible by the definition of an
\hc\ labeling.

Theorem~3.7~\cite{cahit} states that the complete graph $K_n$ is \hc\ for
$n\equiv 0\ (\mod 4)$, and in the proof it is claimed that these are
all possible \hc\ complete graphs. We show in the following theorem
that this is not true.

\begin{thm}
\label{KnHC}
A graph $K_n$ is \hc\ \iff $n\equiv 0,3\ (\mod 4)$, and $n\not=3$.
\end{thm}

\begin{proof}{
If $n\equiv 1,2\ (\mod 4)$ a graph $K_n$ can not be \hc\ by
Lemma~\ref{nechc}. Now if $n\not\equiv 1,2\ (\mod 4)$, we find an
\hc\ labeling for $K_n$.
We know that if $n$ is even, one can decompose $K_n$ into
a $1$--factor and an Eulerian tour.
Now if $n\equiv 0\ (\mod 4)$, the Eulerian tour is \zmc\ and the
$1$--factor is \hc, so $K_n$ is \hc.

Now we consider the case $n\equiv 3\ (\mod 4)$. It is obvious that $K_3$
has an \hc\ labeling. Suppose that $n\ge 7$ and
$u, v, w$ are three vertices of $K_n$. We can find an \hc\ labeling $f$ for
$K_n\setminus\{u,v,w\}$ as in the previous paragraph. The vertices of
$K_n\setminus\{u,v,w\}$ are $n-3$ in number and $(n-3)/2$ of them have
label $+1$. So because $p={(n-3)/4}$ is an integer, we can partition
$V(K_n)\setminus\{u,v,w\}$ into $p$ disjoint subsets $\{a_i,b_i,c_i,d_i\}$,
$i=1,\ldots,p$ such that $f(a_i)=f(b_i)=1$ and $f(c_i)=f(d_i)=-1$,
$\forall i$. We consider two cases to complete the proof.

\begin{figure}[hbt]
\centering
\unitlength 1.00mm
\linethickness{0.4pt}
\begin{picture}(50.58,40.67)
\special{em:linewidth 0.4pt}
\put(4.33,3.67){\circle*{2.50}}
\put(19.33,3.67){\circle*{2.50}}
\put(34.33,3.67){\circle*{2.50}}
\put(49.33,3.67){\circle*{2.50}}
\put(27.00,30.67){\circle*{2.50}}
\put(41.00,38.67){\circle*{2.50}}
\put(11.66,38.67){\circle*{2.50}}
\emline{11.66}{38.67}{1}{41.00}{38.67}{2}
\emline{41.00}{38.67}{3}{49.33}{3.67}{4}
\emline{34.33}{3.67}{5}{41.00}{38.67}{6}
\emline{41.00}{38.67}{7}{4.33}{3.67}{8}
\emline{11.66}{38.67}{9}{27.00}{30.67}{10}
\emline{27.00}{30.67}{11}{19.33}{3.67}{12}
\emline{34.33}{3.67}{13}{27.00}{30.67}{14}
\emline{27.00}{30.67}{15}{49.33}{3.67}{16}
\special{em:linewidth 1.7pt}
\emline{27.00}{30.67}{17}{41.00}{38.67}{18}
\emline{41.00}{38.67}{19}{19.33}{3.67}{20}
\emline{19.33}{3.67}{21}{11.66}{38.67}{22}
\emline{11.66}{38.67}{23}{4.33}{3.67}{24}
\emline{4.33}{3.67}{25}{27.00}{30.67}{26}
\emline{11.66}{38.67}{27}{34.33}{3.67}{28}
\emline{49.33}{3.67}{29}{11.66}{38.67}{30}
\put(2.33,1.34){\makebox(0,0)[cc]{$a_1$}}
\put(17.33,1.67){\makebox(0,0)[cc]{$b_1$}}
\put(32.66,1.67){\makebox(0,0)[cc]{$c_1$}}
\put(47.66,1.67){\makebox(0,0)[cc]{$d_1$}}
\put(9.33,40.67){\makebox(0,0)[cc]{$u$}}
\put(27.33,33.67){\makebox(0,0)[cc]{$v$}}
\put(44.00,40.34){\makebox(0,0)[cc]{$w$}}
\end{picture}
\begin{picture}(50.25,40.67)
\special{em:linewidth 0.4pt}
\put(4.00,3.67){\circle*{2.50}}
\put(19.00,3.67){\circle*{2.50}}
\put(34.00,3.67){\circle*{2.50}}
\put(49.00,3.67){\circle*{2.50}}
\put(11.33,38.67){\circle*{2.50}}
\put(26.67,30.67){\circle*{2.50}}
\put(40.67,38.67){\circle*{2.50}}
\emline{11.33}{38.67}{1}{34.00}{3.67}{2}
\emline{11.33}{38.67}{3}{49.00}{3.67}{4}
\emline{26.67}{30.67}{5}{19.00}{3.67}{6}
\emline{26.67}{30.67}{7}{49.00}{3.67}{8}
\emline{40.67}{38.67}{9}{34.00}{3.67}{10}
\emline{40.67}{38.67}{11}{4.00}{3.67}{12}
\special{em:linewidth 1.7pt}
\emline{11.33}{38.67}{13}{19.00}{3.67}{14}
\emline{11.33}{38.67}{15}{4.00}{3.67}{16}
\emline{26.67}{30.67}{17}{4.00}{3.67}{18}
\emline{26.67}{30.67}{19}{34.00}{3.67}{20}
\emline{40.67}{38.67}{21}{49.00}{3.67}{22}
\emline{40.67}{38.67}{23}{19.00}{3.67}{24}
\put(2.00,1.34){\makebox(0,0)[cc]{$a_i$}}
\put(17.00,1.67){\makebox(0,0)[cc]{$b_i$}}
\put(32.33,1.67){\makebox(0,0)[cc]{$c_i$}}
\put(47.33,1.67){\makebox(0,0)[cc]{$d_i$}}
\put(9.00,40.67){\makebox(0,0)[cc]{$u$}}
\put(27.00,33.67){\makebox(0,0)[cc]{$v$}}
\put(43.67,40.34){\makebox(0,0)[cc]{$w$}}
\end{picture}
\caption{Labeling of $K_n$ when $n\equiv 3\ (\mod 4)$}
\label{KnLbl}
\end{figure}
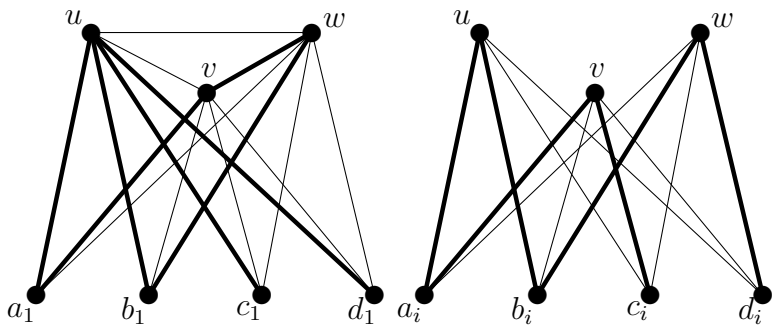

If $p$ is even, we label the un--labeled edges of $K_n$ as in
Figure~\ref{KnLbl}(right) where a thick edge means $(-1)^i$ and a thin one
means $(-1)^{i+1}$; for the edges between $u$, $v$, and $w$ we give to two
of them label $1$ and to the other one label $-1$. If $p$ is odd, for
$i=1$ we use the labels in Figure~\ref{KnLbl}(left) and for $i\ge 2$ we
use the labels in Figure~\ref{KnLbl}(right) where a thick edge means
$(-1)^i$ and a thin one means $(-1)^{i+1}$.
}\end{proof}

\begin{figure}[hbt]
\centering
\special{em:linewidth 0.4pt}
\unitlength .63mm
\linethickness{0.4pt}
\begin{picture}(102.92,57.92)
\special{em:linewidth 0.4pt}
\put(41.67,1.67){\circle*{3.17}}
\put(61.67,1.67){\circle*{3.17}}
\put(61.67,16.67){\circle*{3.17}}
\put(41.67,16.67){\circle*{3.17}}
\put(41.67,21.67){\circle*{3.17}}
\put(41.67,36.67){\circle*{3.17}}
\put(41.67,41.67){\circle*{3.17}}
\put(41.67,56.67){\circle*{3.17}}
\put(61.67,56.67){\circle*{3.17}}
\put(61.67,41.67){\circle*{3.17}}
\put(61.67,36.67){\circle*{3.17}}
\put(61.67,21.67){\circle*{3.17}}
\put(21.17,9.17){\circle*{3.17}}
\put(81.67,9.17){\circle*{3.17}}
\put(81.67,29.17){\circle*{3.17}}
\put(81.67,49.17){\circle*{3.17}}
\put(21.67,49.17){\circle*{3.17}}
\put(21.67,29.17){\circle*{3.17}}
\put(1.67,29.17){\circle*{3.17}}
\put(101.67,29.17){\circle*{3.17}}
\emline{1.67}{29.17}{1}{21.67}{9.17}{2}
\emline{61.67}{16.67}{3}{41.67}{16.67}{4}
\emline{41.67}{16.67}{5}{61.67}{1.67}{6}
\emline{61.67}{1.67}{7}{41.67}{1.67}{8}
\emline{81.67}{9.17}{9}{101.67}{29.17}{10}
\emline{101.67}{29.17}{11}{81.67}{49.17}{12}
\emline{41.67}{56.67}{13}{61.67}{56.67}{14}
\emline{61.67}{56.67}{15}{41.67}{41.67}{16}
\emline{41.67}{41.67}{17}{61.67}{41.67}{18}
\emline{21.67}{49.17}{19}{1.67}{29.17}{20}
\emline{41.67}{21.67}{21}{21.67}{29.17}{22}
\emline{21.67}{29.17}{23}{41.67}{36.67}{24}
\emline{41.67}{36.67}{25}{61.67}{21.67}{26}
\emline{61.67}{21.67}{27}{81.67}{29.17}{28}
\emline{81.67}{29.17}{29}{61.67}{36.67}{30}
\special{em:linewidth 1.7pt}
\emline{41.67}{41.67}{31}{21.67}{49.17}{32}
\emline{21.67}{49.17}{33}{41.67}{56.67}{34}
\emline{41.67}{56.67}{35}{61.67}{41.67}{36}
\emline{61.67}{41.67}{37}{81.67}{49.17}{38}
\emline{81.67}{49.17}{39}{61.67}{56.67}{40}
\emline{41.67}{36.67}{41}{61.67}{36.67}{42}
\emline{61.67}{36.67}{43}{41.67}{21.67}{44}
\emline{41.67}{21.67}{45}{61.67}{21.67}{46}
\emline{41.67}{16.67}{47}{21.67}{9.17}{48}
\emline{21.67}{9.17}{49}{41.67}{1.67}{50}
\emline{41.67}{1.67}{51}{61.67}{16.67}{52}
\emline{61.67}{16.67}{53}{81.67}{9.17}{54}
\emline{81.67}{9.17}{55}{61.67}{1.67}{56}
\emline{81.67}{29.17}{57}{101.67}{29.17}{58}
\emline{1.67}{29.17}{59}{21.67}{29.17}{60}
\end{picture}
\caption{A cubic non--Hamiltonian \hc\ graph}
\label{HCHam}
\end{figure}
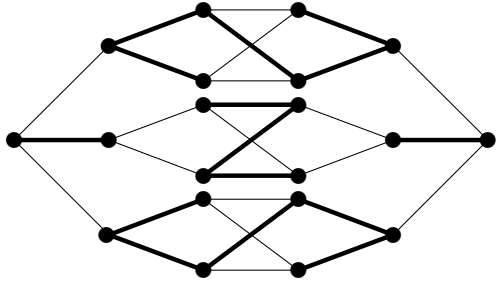

``Theorem~3.9''~\cite{cahit} states that every cubic \hc\ graph is
Hamiltonian
This is not true and the graph shown in Figure~\ref{HCHam} is the
counterexample (for a thin edge we assign the label $-1$ and each thick
edge takes $+1$).

``Theorem~3.10''~\cite{cahit} states that the wheel $W_n$ is \hc\
\iff $n\equiv 1\ (\mod 4)$. In the following theorem we show that
this is not true by giving a necessary and sufficient condition
for \hc ity of a wheel.

\begin{thm}
\label{wheelhc}
The wheel $W_n$ is \hc\ \iff $n$ is odd.
\end{thm}

\begin{proof}{
If $n$ is even, then $W_n$ is not \hc\ by Lemma~\ref{nechc}.
On the other hand if $n$ is odd we give an \hc\ labeling $f$ for $W_n$.
Suppose that $V(W_n)=\{v_0, v_1,\ldots,v_n\}$ and $\Deg v_0=n$.
We define $f(v_0v_i)=f(v_iv_{i+1})=1$ if $1\le i\le n$ and $i$ is even,
also we define $f(v_0v_1)=1$ and for all other edges we give the label~$-1$.
It can easily be seen that $e_f(-1)=e_f(1)=n$, and
$v_f(-1)=v_f(1)=(n+1)/2$. In fact $f(v_i)=1$ for even $i$,
and $f(v_i)=-1$ for odd $i$.
}\end{proof}

\section{Generalizations}

In this section we study another type of graph labeling, called
\hkc\ labeling.

\begin{dfn}
An assignment $f$ of integer labels to the edges of a graph $G$ is called
to be a {\sf\hkc\ labeling}, if for
each edge~$e$ and each vertex~$v$ of $G$ we have $1\le|f(e)|\le k$ and
$1\le|f(v)|\le k$, and for each $i$ with $1\le i\le k$, we have
$|e_f(i)-e_f(-i)|\le 1$ and $|v_f(i)-v_f(-i)|\le 1$.
A graph $G$ is called to be {\sf\hkc}, if it admits a \hkc\ labeling.
\end{dfn}

The following lemma gives a necessary condition for \htc ity of a graph.

\begin{lem}
\label{htcess}
If a graph with an even number of vertices is \htc\ then the number of its
edges is also even.
\end{lem}

\begin{proof}{
If $f$ is a \htc\ labeling for a graph $G$ and $|V(G)|$ is even, then
by lemma~\ref{essential} we have
$2(e_f(1)-e_f(-1)+2e_f(2)-2e_f(-2))=v_f(1)-v_f(-1)+2v_f(2)-2v_f(-2)$.
So $v_f(1)=v_f(-1)$ and since $|V(G)|$ is even, $v_f(2)=v_f(-2)$.
Now $e_f(1)-e_f(-1)+2e_f(2)-2e_f(-2)=0$, so
$e_f(1)-e_f(-1)=e_f(2)-e_f(-2)=0$.
}\end{proof}

\begin{figure}[hbt]
\centering
\scriptsize
\special{em:linewidth 0.4pt}
\unitlength .500mm
\linethickness{0.4pt}
\begin{picture}(65.33,64.83)
\put(62.67,42.33){\circle*{5.00}}
\put(2.67,42.33){\circle*{5.00}}
\put(32.67,62.33){\circle*{5.00}}
\put(32.67,42.33){\circle*{5.00}}
\put(32.67,12.33){\circle*{5.00}}
\put(62.67,2.33){\circle*{5.00}}
\put(2.67,2.33){\circle*{5.00}}
\emline{2.67}{42.33}{1}{2.67}{2.33}{2}
\emline{2.67}{2.33}{3}{32.67}{42.33}{4}
\emline{32.67}{42.33}{5}{32.67}{12.33}{6}
\emline{32.67}{12.33}{7}{62.67}{42.33}{8}
\emline{62.67}{42.33}{9}{62.67}{2.33}{10}
\emline{62.67}{2.33}{11}{32.67}{42.33}{12}
\emline{32.67}{12.33}{13}{2.67}{42.33}{14}
\special{em:linewidth 1.7pt}
\emline{2.67}{42.33}{15}{32.67}{42.33}{16}
\emline{32.67}{42.33}{17}{62.67}{42.33}{18}
\emline{62.67}{42.33}{19}{32.67}{62.33}{20}
\emline{32.67}{62.33}{21}{2.67}{42.33}{22}
\emline{2.67}{2.33}{23}{62.67}{2.33}{24}
\emline{62.67}{2.33}{25}{32.67}{12.33}{26}
\emline{32.67}{12.33}{27}{2.67}{2.33}{28}
\put(5.66,22.00){\makebox(0,0)[cc]{$1$}}
\put(34.00,26.00){\makebox(0,0)[cc]{$2$}}
\put(65.33,21.67){\makebox(0,0)[cc]{$1$}}
\put(10.00,37.33){\makebox(0,0)[cc]{$-1$}}
\put(54.33,37.33){\makebox(0,0)[cc]{$-1$}}
\put(55.33,15.00){\makebox(0,0)[cc]{$-1$}}
\put(9.00,14.33){\makebox(0,0)[cc]{$-1$}}
\put(15.66,10.00){\makebox(0,0)[cc]{$1$}}
\put(46.66,11.00){\makebox(0,0)[cc]{$1$}}
\put(32.00,5.66){\makebox(0,0)[cc]{$1$}}
\put(17.33,45.34){\makebox(0,0)[cc]{$1$}}
\put(46.00,45.67){\makebox(0,0)[cc]{$1$}}
\put(46.00,57.00){\makebox(0,0)[cc]{$1$}}
\put(17.00,56.33){\makebox(0,0)[cc]{$1$}}
\end{picture}
\caption{A counterexample for the converse of Lemma~3.}
\label{htccntr}
\end{figure}
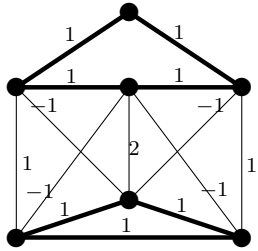

The converse of the above lemma is not necessarily true.
A counterexample is given in Figure~\ref{htccntr}. Note that in place of
bold triangle and quadruple one can put a $C_r$ and a $C_{r+1}$ respectively
for each $r\ge 3$.

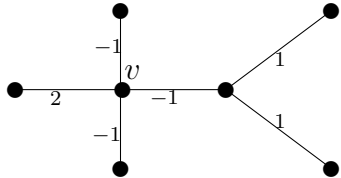
\begin{figure}[hbt]
\centering
\scriptsize
\special{em:linewidth 0.4pt}
\unitlength .70mm
\linethickness{0.4pt}
\begin{picture}(62.58,32.58)
\put(1.33,16.33){\circle*{2.85}}
\put(21.66,16.33){\circle*{2.85}}
\put(41.33,16.33){\circle*{2.85}}
\put(61.33,31.33){\circle*{2.85}}
\put(61.33,1.33){\circle*{2.85}}
\put(21.33,1.33){\circle*{2.85}}
\put(21.33,31.33){\circle*{2.85}}
\emline{1.33}{16.33}{1}{41.33}{16.33}{2}
\emline{41.33}{16.33}{3}{61.33}{31.33}{4}
\emline{61.33}{1.33}{5}{41.33}{16.33}{6}
\emline{21.33}{31.33}{7}{21.33}{1.33}{8}
\put(23.66,19.66){\makebox(0,0)[cc]{\normalsize $v$}}
\put(9.00,14.67){\makebox(0,0)[cc]{$2$}}
\put(19.00,24.33){\makebox(0,0)[cc]{$-1$}}
\put(18.66,7.66){\makebox(0,0)[cc]{$-1$}}
\put(29.66,14.67){\makebox(0,0)[cc]{$-1$}}
\put(51.66,22.00){\makebox(0,0)[cc]{$1$}}
\put(51.66,10.33){\makebox(0,0)[cc]{$1$}}
\end{picture}

\caption{A \htc\ labeling $f$ for which $f^*$ is not \htc}
\label{hk}
\end{figure}

If $f$ is a \hkc\ labeling for a graph $G$, in \cite{cahit} it is
defined another labeling $f^*$ such that $f^*(e)=k+1-f(e)$ if
$f(e)>0$, and $f^*(e)=-k-1-f(e)$ if $f(e)<0$. It is claimed there
that $f^*$ is also a \hkc\ labeling for a graph $G$. There exists
a simple counterexample to the statement. For the tree shown in
Figure~\ref{hk}, the specified labels is a \htc\ labeling, but
$f^*(v)=-5$.

``Theorem~4.2''~\cite{cahit} states that $K_n$ is \htc, \iff
$n\equiv 0\ (\mod 4)$.
We will show this is not true. We prove the following theorem.

\begin{thm}
The complete graph $K_n$ is \htc, if $n\equiv 0,3\ (\mod 4)$, and if
$n\equiv 1\ (\mod 4)$ then $K_n$ is not \htc.
\end{thm}

\begin{proof}{
The \hc\ labelings found in Theorem~\ref{KnHC} are also \htc\ labelings.
On the other hand if $n\equiv 2\ (\mod 4)$ then $K_n$ can not
have a \htc\ labeling by Lemma~\ref{htcess}.
}\end{proof}

The following theorem answers the question of \htc ity of wheels.

\begin{thm}
Every wheel $W_n$ has a \htc\ labeling.
\end{thm}

\begin{proof}{
For odd $n$, we have an \hc\ labeling for $W_n$ by
Theorem~\ref{wheelhc}, which is also an \htc\ labeling. Assume
that $n$ is even and the vertex set of $W_n$ is
$\{v_0,v_1,\ldots,v_n\}$, and $v_0$ is the central vertex (the
vertex with degree $n$). Define $f(v_iv_{i+1})=(-1)^i$ where $1\le
i\le n$ and $v_{n+1}=v_1$. And define $f(v_0v_i)=(-1)^i$ for $2\le
i\le n$, and $f(v_0v_1)=2$. It is straight forward to check that
$f$ is a \htc\ labeling for $W_n$. }\end{proof}

\section*{Acknowledgement}
The authors are very thankful to professor E.\,S.~Mahmoodian for
his useful advices and notes. We also thank Dr.~Ch.~Eslahchi who
read the draft and made useful comments.

\end{document}